\documentclass{article}
\usepackage{amsmath,amssymb,fullpage}
\usepackage{graphicx,subfigure,float,url,color}
\usepackage{mathrsfs}
\usepackage[colorlinks=true]{hyperref}
\usepackage{pdfsync}

\graphicspath{{fig/}}
\def\beq{\begin{equation}}
\def\eeq{\end{equation}}

\def\R{\textrm{I\kern-0.21emR}}
\def\N{\textrm{I\kern-0.21emN}}

\renewcommand{\geq}{\geqslant}
\renewcommand{\leq}{\leqslant}

\newcommand{\vo}{$\dot{V}O2$}
\newcommand{\vom}{$\dot{V}O2_{\mathrm{max}}$}

\newcommand{\Fmax}{F_{\textrm{max}}}

\title{Pace and motor control optimization for a runner}

\author{Amandine Aftalion\footnote{Ecole des Hautes Etudes en Sciences Sociales,  Centre d'Analyse et de Math\'ematique Sociales UMR-8557, Paris, France (\texttt{amandine.aftalion@ehess.fr}).}\and Emmanuel Tr\'elat\footnote{Sorbonne Universit\'e, CNRS, Universit\'e de Paris, Inria, Laboratoire Jacques-Louis Lions (LJLL), F-75005 Paris, France (\texttt{emmanuel.trelat@sorbonne-universite.fr}).}}

\date{}

\begin{document}

\maketitle
\begin{abstract}We present a model which encompasses pace optimization and motor control effort for a runner on a fixed distance. We see that for long races, the long term behaviour is well approximated by a turnpike problem, that allows to define an approximate optimal velocity.  We provide numerical simulations quite consistent with this approximation which leads to a simplified problem. The advantage of this simplified formulation for the velocity is that if we have velocity data of a runner on a race, and have access to his \vom, then we can infer the values of all the physiological parameters. We are also able to estimate the effect of  slopes and ramps.
 \end{abstract}

\section{Introduction}

 The process of running involves a control phenomenon in the human body. Indeed, the optimal pace to run a fixed distance
 requires to use the maximal available propulsive force and energy in order to produce the optimal running strategy.
  This optimal strategy is a combination of cost and benefit: a runner usually wants to finish first or beat the record but
   minimizing his effort. The issue of finding the optimal pacing is a crucial one in sports sciences
   \cite{lapresa,fosterbeating,hettingabrian,hh,hanley2019,thiel2012pacing,tucker2006non,tucker2009physiological} and is still not solved.
    In tactical races, depending on the level of the athlete, and the round on the competition (heating, semi-final or final),
     the strategy is not always the same: the pacing can either be U-shaped (the start and the finish are quicker), J-shaped (greater finishing pace) or reverse J-shaped (greater starting pace) \cite{CasHan,hettingabrian}.

 In this paper, we want  to model this effort minimization as a control problem, solve it and find estimates of the velocity using the turnpike theory of \cite{trelat2020,TZ}.
  We will build on a model introduced by Keller \cite{keller1974optimal}, improved by \cite{aft,AB,AM,AT_RSOS,behncke1993mathematical,bsmall,mathis1989effect,AMH}. The extension by \cite{aft,AB,AM}
   is sufficiently accurate to model real races. We add a motivation equation inspired from the analysis of motor control in the human body \cite{pess}.  This is related to the minimal intervention principle \cite{TJ} so that human effort is minimized through penalty terms.
    We have developed this model for the $200$\,m in \cite{AT_RSOS} and extend it here for middle distance races.

Let us go back to the various approaches based on Newton's second law and energy conservation. Let $d>0$ be the prescribed distance to run. Let $x(t)$ be the position, $v(t)$ the velocity, $e(t)$ the anaerobic energy, $f(t)$ the propulsive force per unit mass. Newton's second law allows to relate force and acceleration through:
 \begin{align*}
& \dot x(t) = v(t)  \qquad\qquad\qquad\qquad\qquad x(0)=0, \qquad x(t_f)=d, \\
& \dot v(t) = -\frac{v(t)}{\tau}+f(t)  \qquad\qquad\qquad v(0)=v^0,
 \end{align*}
 where  $\tau$ is the friction coefficient related to the runner's economy, $t_f$ the final time and $v_0$ the initial velocity. An initial approach by Keller \cite{keller1974optimal} consists in  writing an energy balance: the variation of aerobic energy and  anaerobic energy  is equal to the power developed by the propulsive force, $f(t)v(t)$.
  He assumes that the volume of oxygen per unit of time which is transformed into energy is constant along the race and we call it $\bar \sigma$. If $e^0$ is the initial anaerobic energy, then $\dot e (t)$ is the variation of anaerobic energy and this yields
$$
-\dot e(t) +\bar \sigma =f(t)v(t)\qquad\qquad\qquad e(0)=e^0,\quad e(t)\geq0,\quad e(t_f)=0.
$$
The control problem is to minimize the time $t_f$ to run the prescribed distance $d=\int_0^{t_f} v(t)\ dt$ using a control on the propulsive force $0\leq f(t)\leq f_M$. This model is able to predict times of races but fails to predict the precise velocity profile.

 Experiments have been performed on runners to understand how the aerobic contribution varies with time or distance \cite{hanon2011effects}. Because the available flow of oxygen which transforms into energy  needs some time to increase from its rest value to its maximal value, for short races up to $400$\,m, the function $\sigma$ (which is the energetic equivalent of the oxygen flow) is increasing with time but does not reach its maximal value $\bar \sigma$ or \vom. For longer distances, the maximal value $\bar \sigma$ is reached and $\sigma$ decreases at the end of the race. The longer the race, the longer is the plateau at $\sigma=\bar \sigma$.
  The time when the aerobic energy starts to decrease is assumed to be related to the residual anaerobic supplies \cite{BHKM}. Therefore, in \cite{AB}, to better encompass the link between aerobic and anaerobic effects, the function $\sigma$ is modelled to depend on the anaerobic energy $e(t)$, instead on directly time or distance.
 This leads to the following function $\sigma (e)$ illustrated in Figure \ref{sigma}:\begin{figure}[ht]
\begin{center}
\includegraphics[width=9cm]{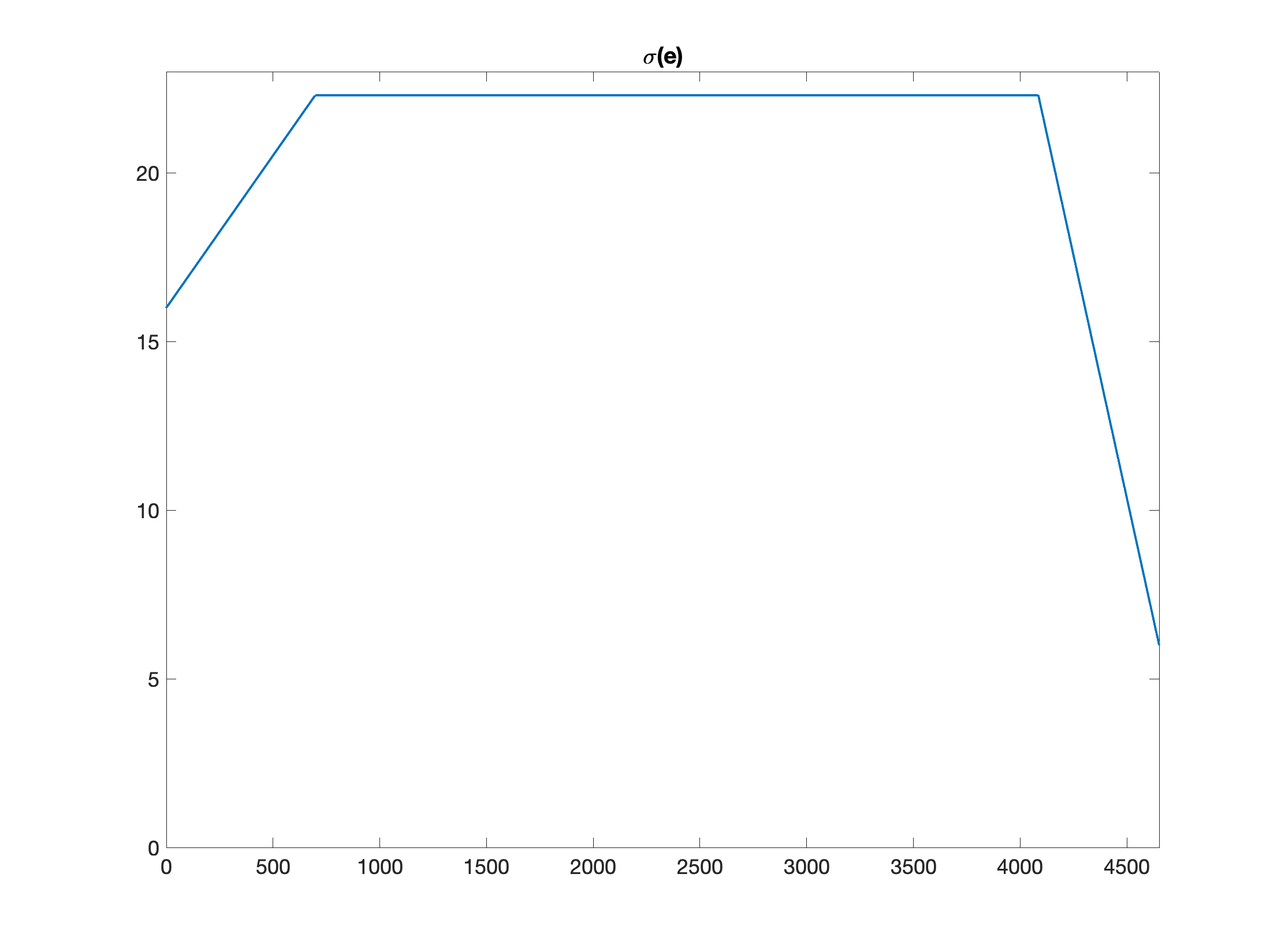}
\end{center}\caption{The function $\sigma (e)$ from \eqref{sigmavar} for $e^0=4651$, $\bar\sigma=22$, $\sigma_f=20$, $\sigma_r=6$, $\gamma_2=566$, $\gamma_1=0.15$.}\label{sigma}\end{figure}
\begin{equation}\label{sigmavar}
\sigma (e)=\left\{\begin{array}{ll}
\displaystyle \bar \sigma \frac{e}{e^0\gamma_1}+\sigma_f \left(1-\frac{e}{e^0 \gamma_1}\right) & \displaystyle\quad \hbox{if}\quad \frac{e}{e^0}<\gamma_1 \\[4mm]
\displaystyle \bar\sigma & \displaystyle \quad \hbox{if}\quad \frac{e}{e^0}\geq \gamma_1 \quad \hbox{and}\quad    e^0 -e\geq \gamma_2\\[2mm]
\displaystyle (\bar \sigma -\sigma_r) \frac{e^0 -e}{\gamma_2}+\sigma_r & \displaystyle \quad\hbox{if}\quad e^0-e<\gamma_2 \end{array}\right.
\end{equation} where $\bar \sigma$ is the maximal value of $\sigma$,  $\sigma_f$ is the final value at the end of the race, $\sigma_r$ is the rest value,
 $e^0$ is the initial value of energy, $\gamma_1 e^0$ is the critical energy at which the rate
 of aerobic energy starts to depend on the residual anaerobic energy and $\gamma_2$ is the energy at which the maximal oxygen uptake $\bar \sigma$ is achieved. Because the anaerobic energy starts at the value $e^0$ and finishes at zero, it depletes in time.
 We observe in our numerical simulations that $e(t)$ decreases, so that $\sigma (e(t))$  and $\sigma (e)$ have opposite monotonicities. The function $\sigma (e(t))$ obtained in our simulations and illustrated in Figure \ref{fig1} is consistent with the measurements of  \cite{hanon2008pacing} or of \cite{hanon2011effects}. The parameters
   $e^0$, $\gamma_1$, $\gamma_2$, $\bar \sigma$, $\sigma_f$, $\sigma_r$ depend on the runner and on the length of the race.

A runner, who speeds up and slows down, chooses to modify his effort. There is a neuro-muscular process controlling human effort. The issue is how to model mathematically this control, coming from motor control or neural drive.
 In Keller's paper \cite{keller1974optimal},  the mathematical control is on the propulsive force. But this yields derivatives of the force which are too big with respect to human ones. Indeed, a human needs some time between the decision to make an effort and the effective change of propulsive force in the muscle. Therefore, in \cite{aft,AB}, the control is the derivative of the propulsive force.
 Nevertheless, putting the control on the derivative of the force seems artificial and it is more satisfactory to actually model the process going from the decision to the muscle. For this purpose, we use the model of
  mechanisms underlying motivation of mental versus physical effort of \cite{pess}. They define the  motor cost of changing a force as the integral of the square of the neural drive $u(t)$. Motor control theory has shown that optimizing this cost minimizes the signal-dependent motor variability and reproduces the cardinal features of movement production. In \cite{pess}, the authors derive the equation for the derivative of the force which limits the variation of the force through the neural drive $u(t)$:
   \begin{itemize}
   \item the force increases with the neural drive so that $\dot f$ is proportional to $u$;
   \item the force is bounded by a maximal force even when the neural drive increases so that
   $\dot f$ is proportional to $u(\Fmax -f)$;
   \item without excitation, it decreases exponentially so that $\dot f$ is proportional to $u(\Fmax -f)-f$;
   \item the dynamics of contraction and excitation depends on the muscular efficiency $\gamma$ so that $\dot f$ is proportional to $\gamma$.
   \end{itemize}
Therefore,  following \cite{pess}, and as in \cite{AT_RSOS}, we add an equation for the variation of the force.
This leads to the following system:
\begin{align}\label{equ}
& \dot x(t) = v(t)  \qquad\qquad\qquad\qquad\qquad\qquad\qquad\quad x(0)=0, \quad x(t_f)=d, \\
& \dot v(t) = -\frac{v(t)}{\tau}+f(t)\label{equv}
 \qquad\qquad\qquad\qquad\qquad\quad v(0)=v^0, \\
& \dot f(t) = \gamma \Big( u(t) (\Fmax-f(t)) - f(t) \Big) \qquad\qquad  f(t)\geq 0,  
\label{equf}\\
& \dot e(t) = \sigma(e(t)) -f(t)v(t)\qquad\qquad\qquad\qquad\quad e(0)=e^0,\quad e(t)\geq0,\quad e(t_f)=0,\label{eqe}
\end{align}
where $e^0$ is the initial energy, $\tau$ the friction coefficient related to the runner's economy, $F_{\mathrm{max}}$ is a threshold upper bound for the force, $\gamma$ the time constant of motor activation and $u(t)$ the neural drive which will be our control.
We observe in our simulations that, in order to minimize the time, the force $f(t)$ remains positive along the race without the need to put it as a constraint.
Let us point out that it follows from Equation \eqref{equf} that $f(t)$ cannot cross $F_{\mathrm{max}}$ increasing. Therefore, with our choice of parameters (the value of $e^0$ is not large enough), we observe that $f(t)$ always remains below $F_{\mathrm{max}}$ without putting any bound on the maximal force.
In this paper,  we do not take into account the effect of bends because for long races, they have minor effects on the velocity.

The optimization problem consists in minimizing the difference between the cost and the benefit. In \cite{pess}, the expected  cost is proportional to the motor control which is the $L^2$ norm of the neural drive $u(t)$.
  On the other hand, the benefit is proportional to the reward, and can be estimated for instance to be proportional to $-t_f$. Indeed, one could imagine the reward is a fixed amount to which is subtracted a number proportional to the difference between the world record and the final time.
 Similarly, one could add other benefits or costs linked to multiple attempts or the presence of a supporting audience. One could think of adding other costs, for instance in walking modeling, the cost is proportional to the jerk, which is the $L^2$ norm of the derivative of the centrifugal acceleration \cite{laumond,capo}.
 In this paper, we choose to model the simplest case where the benefit is the final time and the cost is the motor control. This leads to the following minimization:
 \begin{equation}\label{optcond} \min \left( t_f+\frac{\alpha}{2} \int_0^{t_f} u(t)^2\, dt \right)\end{equation} where $\alpha>0$ is a weight to be determined so that the second term is a small perturbation of the first one, and therefore both terms are minimized.

As soon as the race is sufficiently long (above $1500$\,m), one notices (see \cite{hanon2011effects} and our numerical simulations) the existence of a limiting problem where $v$ and $f$ are constant and $e$ is linearly decreasing. Therefore, it is natural to expect
that the turnpike theory of \cite{TZ} (see also \cite{trelat2020}) provides very accurate estimates for the mean velocity, force and the energy decrease. The turnpike theory in optimal control stipulates that, under general assumptions, the optimal solution of an optimal control problem in sufficiently large fixed final time
 remains essentially constant, except at the beginning and at the end of the time-frame. We refer the reader to \cite{TZ} for a complete state-of-the-art and bibliography on
the turnpike theory. Actually, according to \cite{TZ}, due to
  the particular symplectic structure of the first-order optimality system derived from the Pontryagin maximum principle, the optimal state, co-state (or adjoint vector) and optimal control   are, except around the terminal points, exponentially close to steady-states, which are themselves the optimal solutions of an associated \emph{static} optimal control problem. In this result, the turnpike set is a singleton, consisting of this optimal steady-state which is of course an equilibrium of the control system. This is the so-called \emph{turnpike phenomenon}.
A more general version has recently been derived in \cite{trelat2020}, allowing for more general turnpike sets and establishing a turnpike result for optimal control problems in which some of the coordinates evolve in a monotone way while some others are partial steady-states.
This result applies to our problem and we want to use it to simplify the runner's model for potential software applications.

 The paper is organized as follows. Firstly, we present numerical simulations of \eqref{equ}-\eqref{equv}-\eqref{equf}-\eqref{eqe}-\eqref{optcond}, then we describe our simplified problem and how to derive it. In Section 4, we study a more realistic \vo\ and in Section 5, the effects of slopes.

\section{Numerical simulations}\label{sec2}
Optimization and numerical implementation of the optimal control problem \eqref{equ}-\eqref{equv}-\eqref{equf}-\eqref{eqe}-\eqref{optcond} are done by combining automatic differentiation softwares with the modeling language AMPL~\cite{Fourer2002} and expert optimization routines with the open-source package IpOpt~\cite{Waechter2006}.
This allows to solve for the velocity $v$, force $f$, energy $e$ in terms of the distance providing the optimal strategy and the final time.
 As advised in \cite{trelat2020,TZ}, we initialize the optimization algorithm at the turnpike solution that we describe below.

We have chosen numerical parameters to match the real race of $1500$\,m described in \cite{hanon2008pacing} so that $d=1500$. The final experimental time for real runners is $245$\,s.
 The runners are middle distance runners successful in French regional races. Their \vom\ is around $66$\,ml/mn/kg.  Because it is estimated that one liter of oxygen produces an energy of about $21.1$\,kJ via aerobic cellular mechanisms \cite{Per}, the energetic equivalent of $66$\,ml/mn/kg is $66\times 21.1$\,kJ/mn/kg. Since we need to express $\sigma$, the energetic equivalent of \vo\ in SI units, we have to turn the minutes into seconds and this provides an estimate of the available energy per $kg$ per second which is $ 66/60\times 21.1 \simeq 22$.
 This leads to a maximum value $\bar\sigma=22$ of $\sigma$.    From \cite{hanon2008pacing}, the decrease in \vo\ at the end of the race is of about $10\%$ when the anaerobic energy left is $15\%$.
  Therefore, we choose  the final value of $\sigma$ to be $10\%$ less than the maximal value, that is $\sigma_f=20$, and  $\gamma_1=0.15$. To match the usual rest value of \vo, we set   $\sigma_r = 6$. The other parameters are identified so that the solution  of \eqref{equ}-\eqref{equv}-\eqref{equf}-\eqref{eqe}-\eqref{optcond} matches the velocity data of \cite{hanon2008pacing}:
 $\gamma_2 = 566$, $\alpha=10^{-5}$,
 $F_{\mathrm{max}} = 8$, $\tau= 0.932$, $ e^0 = 4651$, $\gamma=0.0025$, $ v^0 = 3$. Let us point out that our model of effort is not appropriate to describe the very first seconds of the race. Therefore, we choose artificially $v_0=3$ which allows, with our equations, to have a more realistic curve for the very few points, than starting from $v_0=0$. Otherwise, one would need to refine the model for the start.

  In \cite{pess}, the equivalent of $\alpha$ is determined by experimental data.
 In our case, we  have noticed that, depending on $\alpha$, either $u$ is negative with a minimum or changes sign with a minimum and a maximum. Also, when $\alpha$ gets too small,  $\dot f$ is almost constant.
 The choice of $\alpha$ is made such that the second term of the objective is a small perturbation of the first one, and can act at most on the tenth of second for the final time.

\begin{figure}[ht]
\centerline{\includegraphics[width=15.1cm]{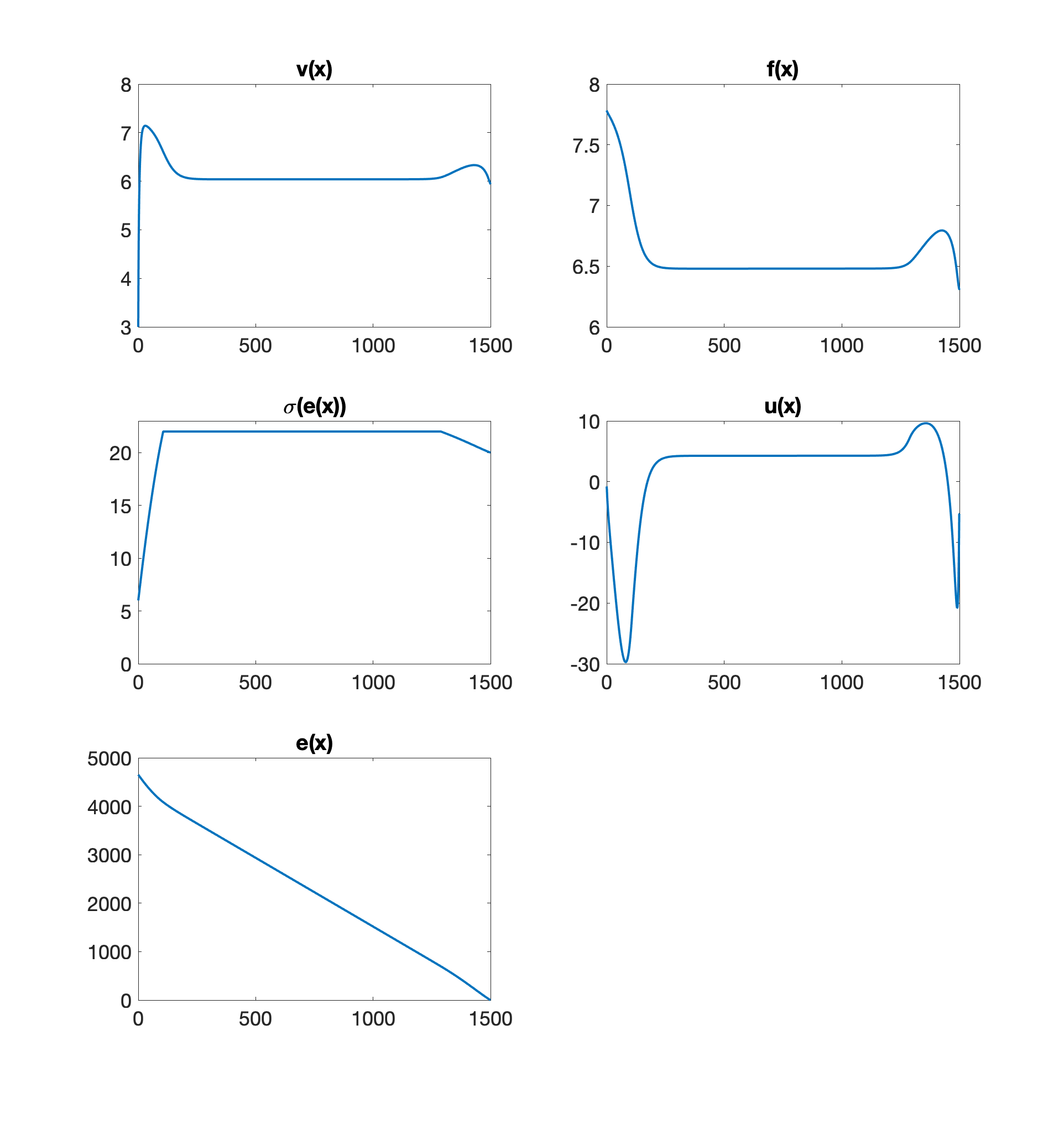}}
\caption{Velocity $v$, force $f$,  energetic equivalent of the oxygen uptake $\sigma(e)$, motor control $u$ and energy $e$ vs distance on a $1500$\,m. All functions (except $e$) display a plateau in the middle of the race corresponding to the turnpike phenomenon, except the energy which is affine. In this numerical simulation, the duration of the race is $244$\,s.}\label{fig1}
\end{figure}

  With these parameters, we simulate the optimal control problem \eqref{equ}-\eqref{equv}-\eqref{equf}-\eqref{eqe}-\eqref{optcond} and plot the velocity $v$, the propulsive force $f$, the motor control $u$, the energetic equivalent of the oxygen uptake $\sigma(e)$,  and the anaerobic energy $e$ vs  distance in Figure \ref{fig1}. Though they are computed as a function of time, we find it easier to visualize them as a function of distance.
 The velocity increases until reaching a peak value, then decreases to a mean value, before the final sprint at the end of the race. This is consistent with usual tactics which consist in an even pace until the last $300$\,m  where the final sprint starts. This final sprint takes place when the function $\sigma(e(t))$ starts decreasing. The function $\sigma$ is the energetic equivalent of \vo. It increases to its plateau value, then decreases at the end of the race when the anaerobic supply gets too low. The control $u$ also has a plateau at the middle of the race leading to a plateau for the force as well. The velocity and force follow the same profile. The energy is decreasing and almost linear when the velocity and force are almost constant.

In Figure \ref{fig1}, we point out that we obtain  an almost steady-state in the central part of the race for the motor control, the force and the velocity. We find from Figure \ref{fig1} the central value for the motor control $ u_{\textrm{turn}}=4.26$, the force $f_{\textrm{turn}}=6.48$ and the velocity $v_{\textrm{turn}} = 6.04$.  We want to analyze this limit analytically. We will also try to construct local models for the beginning and end of the race.

\section{Main results using turnpike estimates}

The optimal control problem \eqref{equ}-\eqref{equv}-\eqref{equf}-\eqref{eqe}-\eqref{optcond}  involves a   state variable, namely, the energy $e(t)$, which goes from $e^0$ to 0, and thus has no equilibrium. The turnpike theory has been extended in \cite{trelat2020} to this situation when the steady-state is replaced by a partial steady-state (namely, $v$ and $f$ are steady),
  and $e(t)$ is approximated by an affine function satisfying the imposed constraints $e^0$ at initial time and $0$ at final time. In what follows,
   we denote the approximating turnpike trajectory with an upper bar, corresponding to a constant function $\sigma (e)=\bar \sigma$.
More precisely,  we denote by $t\mapsto(\bar v_c, \bar e_c(t), \bar f_c)$ the turnpike trajectory defined on the interval $[0,\bar t_c]$ so that $\bar v_c$
 and $\bar f_c$ are steady-states (equilibrium of the control dynamics \eqref{equ}-\eqref{equv}-\eqref{equf}) with $\bar v_c
=\bar f_c \tau$, and $\bar e_c(t)$ is affine: $$ \dot{\bar e}_c (t)=\bar\sigma -\bar f_c\bar v_c$$ and satisfies the terminal constraints
  $\bar e_c(0)=e^0$ and $\bar e_c(\bar t_c)=0$, while $d=\bar v_c \bar t_c$. Integrating yields
\begin{equation}\frac{\bar v_c^2}\tau -\bar \sigma =e^0
\frac {\bar v_c}d.\label{barv} \end{equation}
The mean velocity $\bar v_c$ can be solved from \eqref{barv} to get
\begin{equation}\bar v_c =\frac{e^0 \tau}{2d} +\sqrt{\bar\sigma\tau+\left ( \frac{e^0\tau}{2d}\right )^2}.\label{eqvbar}\end{equation} We observe
 that the value of $\bar v_c$ increases with $e^0$, $\tau$ (which is the inverse of friction) and $\bar \sigma$, but is not related to the maximal force. Indeed, the maximal propulsive force controls the acceleration at the beginning and end of the race, but not the mean velocity in the middle of the race.
In the case of our simulations, $\bar v_c= 6.2$ which is slightly overestimated with respect to the simulation value $v_{\textrm{turn}}=6.04$.

\medskip

    We next elaborate to show how the turnpike theory can be applied to the central part of the race where $\sigma$ is constant and allows to derive very accurate approximate solutions.

 If one takes into account the full shape of $\sigma(e)$, made up of three parts, then  the velocity curve is made up of three parts. In the rest of the paper, we will derive the following approximation for the velocity:

\medskip
\centerline{
\boxed{\begin{minipage}{0.99\textwidth}
\begin{equation}\label{velocity}
v(t)=\left\{\begin{array}{ll}
\displaystyle v_0e^{-t/\tau}+\left (v_{\mathrm{max}}+\frac t{t_1}(\bar v -v_{\mathrm{max}})\right )(1-e^{-t/\tau}) & \quad\hbox{if}\quad 0\leq t \leq t_1, \\[2mm]
\displaystyle \bar v &\quad\hbox{if}\quad t_1\leq t \leq t_2,\\[2mm]
\displaystyle \frac{ \tau \Fmax}{ 1+  ( {\Fmax}/ {\bar f}-1) e^{-\gamma \lambda  \Fmax (t-t_2)}}
&\quad\hbox{if}\quad t_2\leq t \leq t_f.
 \end{array}\right.
\end{equation}\end{minipage}} }
\bigskip

\noindent
The parameters appearing in the formula are defined as follows: $v_0$ is the initial velocity in \eqref{equv}, $\bar v$ is obtained as the positive root that is bigger than $\sqrt{\bar \sigma \tau}$ of
 \begin{equation}\label{dv}
d=\frac{\bar v\gamma_2}{\frac {\bar v^2}\tau-\sigma_r}+\bar v\frac{e^0(1-\gamma_1)-\gamma_2}{\frac {\bar v^2}\tau-\bar\sigma}+\frac{\bar ve^0\gamma_1}{\frac {\bar v^2}\tau-\sigma_f},
\end{equation}
 $t_1$ is given by
\begin{equation}\label{tt1}
t_1=\frac{\gamma_2} {\frac{\bar v^2}\tau -\sigma_r},
\end{equation}
$v_{\max}=f^0\tau$, where $f^0$ is the positive root of the trinomial
\begin{multline}\label{f00}
\int_0^{t_1}  \left(f^0+t\frac{\bar v/ \tau-f^0}{t_1}\right)\left(v_0 e^{-t/\tau}+\left(\tau f^0+t\frac{\bar v-\tau f^0}{t_1}\right)(1-e^{-t/\tau})\right )e^{\frac {\bar \sigma -\sigma_r} {\gamma_2}(t-t_1)} \, dt \\
= \gamma_2+\frac{\sigma_r\gamma_2}{\bar \sigma -\sigma_r}\left(1-e^{-\frac {\bar \sigma -\sigma_r} {\gamma_2}t_1}\right)  ;
\end{multline}
from this, we compute $d_1=\int_0^{t_1} v(t) \, dt$. We define $\bar d=\frac{e^0(1-\gamma_1)-\gamma_2}{\frac {\bar v^2}\tau-\bar\sigma}$, the length of the turnpike, and
\beq\label{dtend} \Delta t_{\textrm{end}}= \frac {d-d_1-\bar d}{\bar v};\eeq
   $\lambda$ is chosen such that, if $A=\frac {\bar \sigma -\sigma_r}{\gamma_1 e^0}$, then there is an $L^2$ estimate for the velocity at the end of the race:
\beq \label{eqenerv}
\int_{0}^{\Delta t_{\textrm{end}}} \left(\frac{ \tau \Fmax}{ (1+  ( {\Fmax}/ {\bar f}-1) e^{-\gamma \lambda  \Fmax t})} \right)^2 e^{-At}\ dt=\tau \frac{\sigma_f}A (1-e^{-A \Delta t_{\textrm{end}}})+
\tau \gamma_1 e^0; \eeq
moreover, the time $t_2$ is defined so that
\beq\label{t2t1}
t_2-t_1=\frac{1}{\bar v} \left( d-\int_0^{t_1}v(t)\ dt -\int_0^{\Delta t_{\textrm{end}}}\frac{\tau \Fmax}{ 1+  ( {\Fmax}/ {\bar f}-1) e^{-\gamma \lambda  \Fmax t}}\ dt\right) .
\eeq
and $t_f=t_2+\Delta t_{\textrm{end}}$.

\medskip

Let us explain the general meaning of these computations. Equation \eqref{dv} is based on the hypothesis that $v$ and $f$ are constant values and uses the shape of $\sigma$ and the energy equation to compute the duration and length of each phase.  From the first phase, we derive the value of $t_1$ in \eqref{tt1}. Then we compute the initial force that corresponds to the correct energy expenditure in the first phase through \eqref{f00}. This provides, through the integral of the velocity the distance $d_1$ of the first phase. We next approximate the distance and time of the last phase using the distance and time of turnpike through \eqref{dtend}. Once we have the duration of the last phase, we again match the energy expenditure in \eqref{eqenerv}. This provides the velocity profile of the last phase and therefore the distance of the last phase. In order to match the total distance, we have to slightly modify the length of the central turnpike part in \eqref{t2t1}. From the computational viewpoint, these steps correspond to the first successive approximations in the Newton-like solving of a system of nonlinear equations.

\medskip

The velocity curve \eqref{velocity} goes from the initial velocity $v_0$ to a maximum velocity, then down to $\bar v$, which is the turnpike value. At the end of the race, the velocity increases to the final velocity.
 This type of curve is quite consistent with velocity curves in the sports literature, see for instance \cite{fosterbeating,hanley2019}, and with our simulations illustrated in Figure \ref{fig1}.

We see that $t_1$ increases with $\gamma_2$, while
$t_f-t_2$ increases with $\gamma_1$.

For the values of parameters of Section \ref{sec2}, we find from \eqref{dv} that $\bar v=6.06$, which is to be compared to the value in Figure \ref{fig1}, $v_{\textrm{turn}}=6.04$.
 Then from \eqref{tt1} $t_1=16.95$, from \eqref{f00} that $f^0=8.2$, $d_1=111.84$. We deduce from \eqref{dtend} $\Delta t_{\textrm{end}}=34.42$, from \eqref{eqenerv}, $\lambda=0.39$, from \eqref {t2t1} $t_2=210.76$, $t_f=245.19$ (very close to the $244$\,s obtained in the numerical simulation in Figure \ref{fig1} and to the experimental value of $245$\,s) and we find $v_f=6.33$ at the final time. We point out that in the turnpike region, this yields $\bar f=\bar v/\tau=6.5$ and $\bar u=\bar f/(\Fmax -f)= 4.34$, very close to the values in Figure \ref{fig1}, $f_{\textrm{turn}}=6.48$ and $ u_{\textrm{turn}}=4.26$.

\medskip

\begin{figure}[ht]
\centerline{\includegraphics[width=18.1cm]{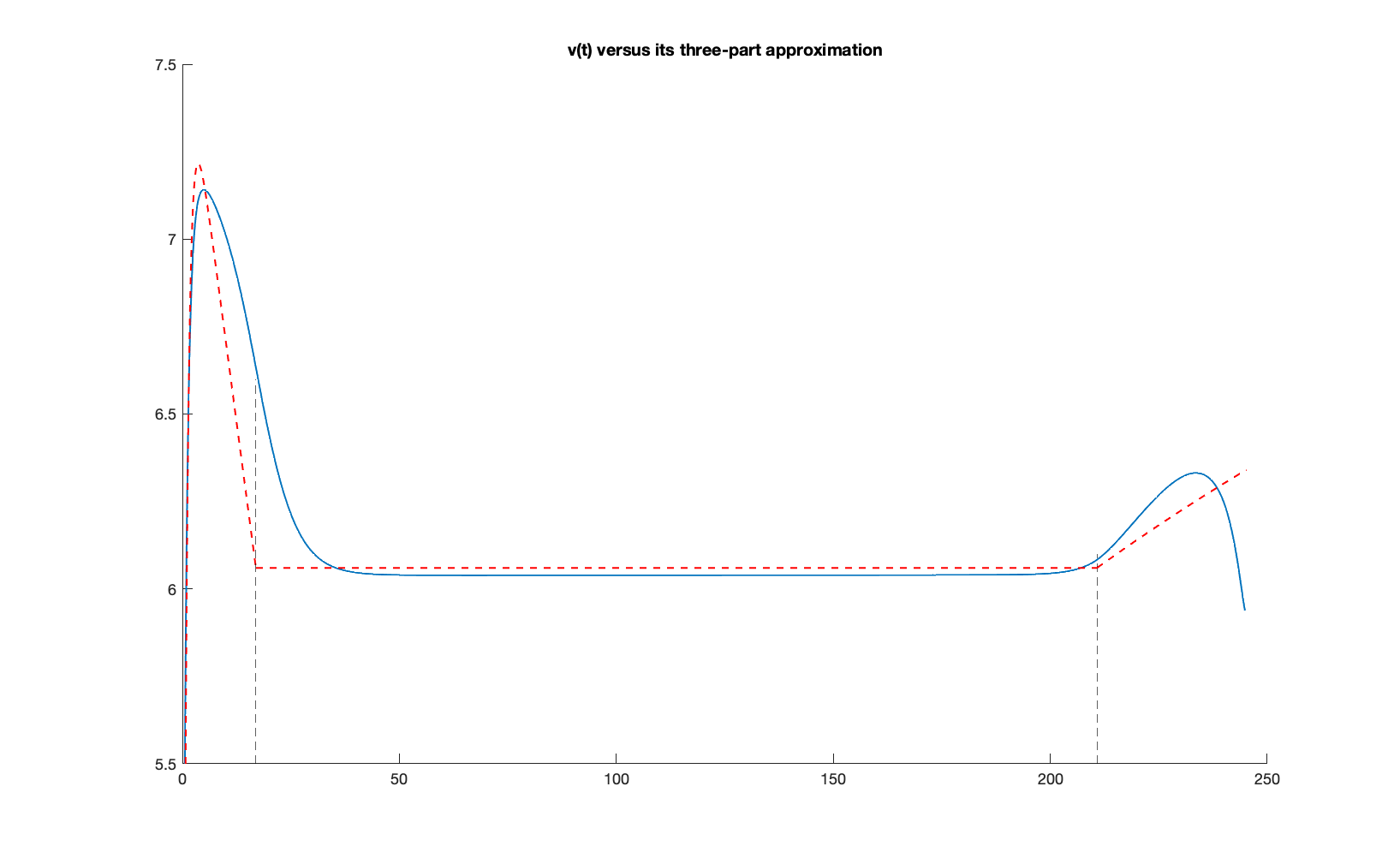}}
\caption{Velocity $v$ as a solution of the simulation (blue) of \eqref{equ}-\eqref{equv}-\eqref{equf}-\eqref{eqe}-\eqref{optcond} and approximate solution given by \eqref{velocity} (red).}\label{figvapp}
\end{figure}

We have illustrated in Figure \ref{figvapp} the approximate solution \eqref{velocity} together with the numerical solution of the full optimal control problem \eqref{equ}-\eqref{equv}-\eqref{equf}-\eqref{eqe}-\eqref{optcond}. We see that the duration of the initial phase is slightly underestimated,
while the duration of the final phase is very good. The estimate of the sprint velocity at the end is also very good.
Note that the simulation of the full optimal control problem produces a decrease of velocity at the very end of the race which is not captured by our approximation, but this changes very slightly the estimate on $t_f-t_2$ or on the sprint velocity at the end and is not meaningful for a runner, so we can safely ignore it for our approximations.

\medskip

The advantage of formulation \eqref{velocity} is that if we have velocity data of a runner on a race, and have access to his \vom, that is $\bar \sigma$, then we can infer the values of all the physiological parameters: from the velocity curve at the beginning, we can determine $\tau$ and $v_{\max}$. The value of $\bar v$ and \eqref{eqvbar} yield $e^0$. From the values of $t_1$ and $t_2$, we deduce $\gamma_1$ and $\gamma_2$. In order to have more precise values, we can always perform an identification of the parameters using the full numerical code, but from these approximate values, we have enough information to determine the runner's optimal strategy on other distances.

The rest of the section is devoted to deriving \eqref{velocity}.

\hfill


\subsection{Central turnpike estimate}
In the central part of the race, $\sigma(e)=\bar \sigma$ is constant. Therefore in this part, when $e(t)$ is between $e^0-\gamma_2$ and $\gamma_1 e^0$, we can apply the turnpike theory of \cite{trelat2020}.
Then we have $v(t)\simeq\bar v$, $f(t)\simeq\bar f$, $u(t)\simeq\bar u$ with
$$
\bar f=\frac {\bar v}\tau,\quad \bar u = \frac{\bar f}{F_{\mathrm{max}}-\bar f}.
$$
We have to integrate
$$
\dot{\bar e}(t)=\sigma(\bar e(t))-\frac {\bar v^2}\tau,\qquad \bar e(t_1)=e^0-\gamma_2\quad \bar e(t_2)=\gamma_1 e^0.
$$
 We find
$$
e^0(1-\gamma_1)-\gamma_2=(t_2- t_1) \left({\frac {\bar v^2}\tau-\bar \sigma}\right) .
$$ This is consistent with \eqref{barv} which is the same computation but on the whole interval, that is with $\gamma_1=\gamma_2=0$. The value for $t_2-t_1$ is $194.64$.

As a first approximation, we can assume that on the two extreme parts of the race, $v$ and $f$ can be taken to be constants. We will see below why this assumption is reasonable. Therefore we can solve
$$
\dot{\bar e}(t)=\sigma(\bar e(t))-\frac {\bar v^2}\tau \qquad\qquad \bar e(0)=e^0,\quad \bar e(t_1)=e^0-\gamma_2,
\quad \bar e(t_2)=\gamma_1 e^0,\quad e(\bar t)=0.$$  Therefore, $\bar t$ is the final time of the turnpike trajectory  defined by $\bar e(\bar t)=0$.
The initial and final parts of the race produce exponential terms, namely
\begin{equation}\label{tint}
\frac {\bar \sigma -\sigma_r}{\frac {\bar v^2}\tau-\sigma_r}=1-e^{-\frac {(\bar \sigma-\sigma_r)t_1}{\gamma_2}}
\qquad\quad\hbox{and}\qquad\quad
\frac {\bar \sigma -\sigma_f}{\frac {\bar v^2}\tau-\sigma_f}=1-e^{-\frac {(\bar \sigma-\sigma_f)(\bar t -t_2)}{e^0\gamma_1}}  .
\end{equation}
Therefore, for the total distance $d$, we find, summing our estimates,
\begin{equation}\label{tbarnew}
\bar t=\frac d{\bar v}=\frac{e^0(1-\gamma_1)-\gamma_2}{\frac {\bar v^2}\tau-\bar\sigma}-\frac{\gamma_2}{\bar\sigma -\sigma_r} \ln \left ( 1-\frac {\bar \sigma -\sigma_r}{\frac {\bar v^2}\tau-\sigma_r}\right)-\frac {e^0\gamma_1}{\bar\sigma -\sigma_f} \ln \left ( 1-\frac {\bar \sigma -\sigma_f}{\frac {\bar v^2}\tau-\sigma_f}\right).
\end{equation}
If the initial and final parts are not too long, then \eqref{tint}  can be approximated by
\begin{equation}\label{t1t2}
t_1\simeq \frac{\gamma_2}{\frac {\bar v^2}\tau-\sigma_r}
\qquad\quad\hbox{and}\qquad\quad
\bar t -t_2\simeq \frac{e^0\gamma_1}{\frac {\bar v^2}\tau-\sigma_f}
\end{equation}
and therefore, from \eqref{tbarnew}, $\bar v$ can be approximated by \eqref{dv}.
For the values of parameters of Section \ref{sec2}, \eqref{dv} yields $\bar v=6.06$.  The intermediate times can be computed from \eqref{t1t2}: $t_2=35.96$\,s and  $t_1=16.95$\,s.
 This also yields the distances of each part by multiplying by $\bar v$. In the following, we will keep this value of $t_1$ but improve the estimate for $t_2$.

\hfill

Note that this turnpike calculation can be used the other way round: if one knows the mean velocity, $d$, $\tau$ and $\bar \sigma$, it yields an estimate of the energy $e^0$ used while running, as well as the aerobic part which is $\bar \sigma d/\bar v$.

The next step is to identify reduced problems for the beginning (interval $(0,t_1)$) and end of the race (interval $(t_2,\bar t )$). The two are not totally equivalent since at the beginning we have an initial condition for the velocity $v$  whereas on the final part the final velocity is free.

\subsection{Estimates for the beginning of the race}

The  problem is to approximate  the equations for $v$, $f$, $e$  with boundary conditions
$$
v(0)=v^0,\quad v(t_1)=\bar v,\quad f(t_1)=\bar f,\quad e(0)=e^0,\quad e(t_1)= e^0-\gamma_2.
$$
Here, $f(0)$ is free.

We integrate the energy equation and find
$$\int_0^{t_1} f(t) v(t) \ dt=\int_0^{t_1} \left(\dot e(t) -\sigma (e(t))\right) dt.$$
In this regime, $\sigma(e)$ is linear, and this equation can be integrated explicitly. Indeed, let $A=\frac{\bar \sigma -\sigma_r}{\gamma_2}$, then
\beq\label{enereq}
-\gamma_2=\frac{\sigma_r\gamma_2}{\bar \sigma -\sigma_r}(1-e^{-At_1})-e^{-At_1}\int_0^{t_1}f(t)v(t)e^{At}\ dt.\eeq
Because we are in a regime of parameters where $At$ is small, we can expand the exponential terms. The approximation which consists in assuming that the integral of $fv$ can be approximated by the mean value of $fv$ is good, and therefore this justifies the turnpike estimate of the previous section and this yields the estimate \eqref{tt1} of $t_1$.

  Now let us assume $t_1$ is prescribed.
If we fix the interval $(0,t_1)$, we have the equations for $v$ and $f$ with
\beq\label{initbeg2}v(0)=v^0,\quad v(t_1)=\bar v,\quad f(0)=f^0,\quad f(t_1)=\bar f.\eeq
Here $f^0$ is unknown and we want to minimize the motor control only. For this part, we can assume that the minimization of the motor control leads to a linear function $f$ as explained in the Appendix.
Therefore, $f(t)=f^0+t(\bar f-f^0)/t_1$ and  $v(t)=v_0 e^{-t/\tau}+\tau f(t)(1-e^{-t/\tau})$ to approximate \eqref{equv}.
We plug this into \eqref{enereq}  and then we find that $f^0$  is a solution of \eqref{f00}.
This can be integrated analytically or numerically to determine $f^0$. In our case, $f^0=8.2$. This yields the first line of \eqref{velocity}  with $v_{\mathrm{max}}=\tau f^0$.

\subsection{End of the race}
 Once the beginning and central part of the race are determined, the duration of the end of the race is determined so that the prescribed distance $d$ is run through \eqref{dtend}.

\begin{figure}[ht]
\centerline{\includegraphics[width=18cm]{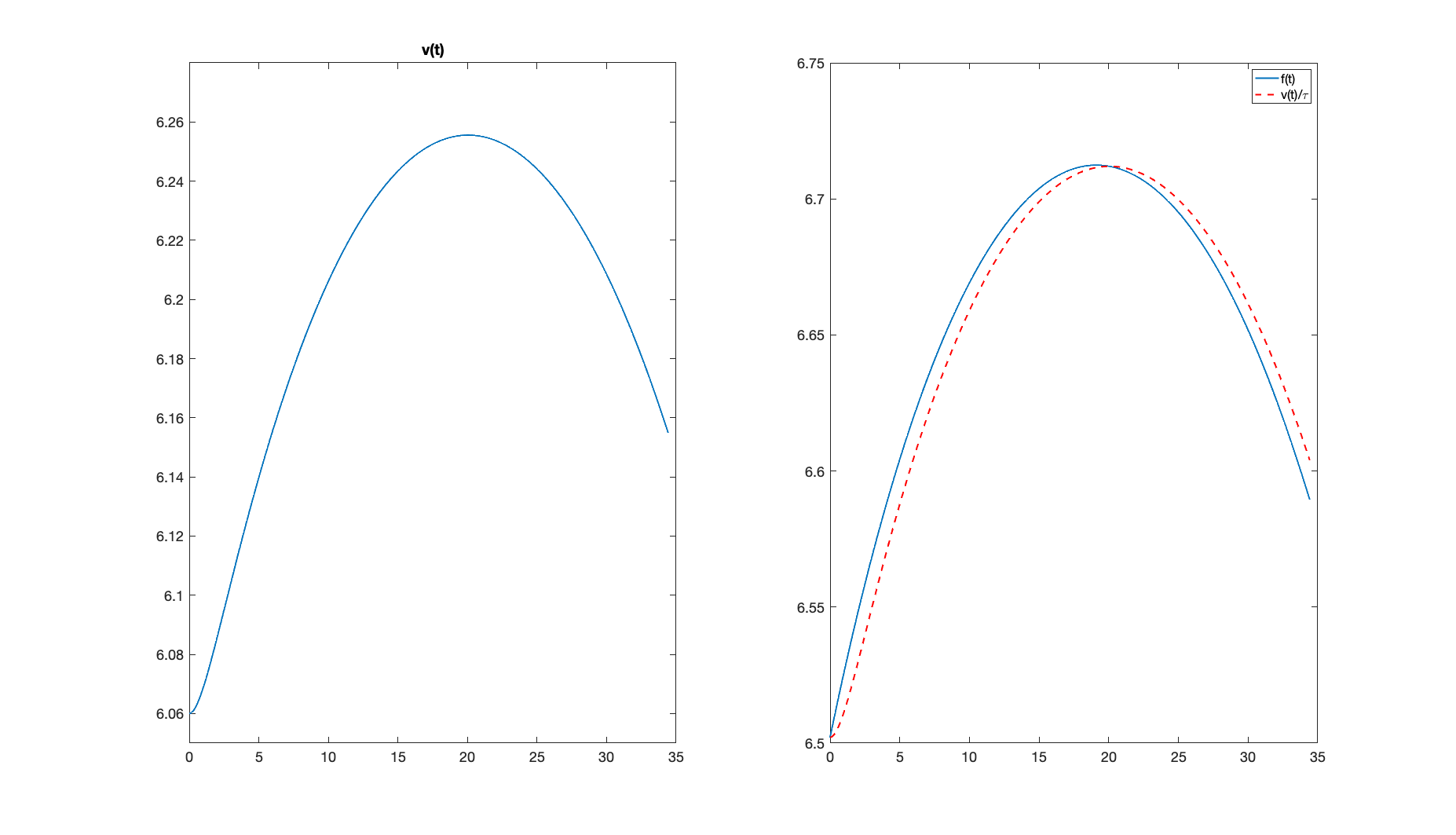}}
\caption{Velocity and force solving the equations for $v$, $f$, $e$ on the interval $(t_2,t_f)$ with initial and final values
\eqref{bcend}. The force $f(t)$ is compared to the value $v(t)/\tau$.}\label{figfin}
\end{figure}

The problem describing the end of the race consists in solving the equations for $v$, $f$, $e$ on the interval $(t_2,t_f)$ with initial and final values
\beq\label{bcend} v(t_2)=\bar v,\quad f(t_2)=\bar f,\quad e(t_2)=\gamma_1 e^0,\quad e(t_f)=0.\eeq
  This yields the simulation in Figure \ref{figfin}. We observe that $f(t)$ and $v(t)/\tau$ are very close, as expected.

In the following, we will assume  that $\dot v$ is negligible in front of $v/\tau$, so that $v\simeq f\tau$, which removes an equation. Then using the specific shape of $\sigma$, the energy equation becomes, denoting $A=\frac{\bar\sigma-\sigma_f}{e^0\gamma_1}\simeq 0.0028$,
$$
e^{A(t-t_2)}\frac{d}{dt}{{\left(e(t)e^{-A(t-t_2)}\right)}}=\sigma_f-\tau f(t)^2.
$$
Then we need to integrate this energy equation and find
\beq \label{eqenerA} \tau \int_{t_2}^{t_f} f(t)^2 e^{-A(t-t_2)}\, dt=\frac{\sigma_f}A (1-e^{-A(t_f-t_2)})+
\gamma_1 e^0.\eeq
%
 The reduced optimal control problem for the end of the race is therefore
 \beq \label{eqfend}
 \begin{split}
 & \min \int_{t_2}^{t_f} u(t)^2\, dt \\
 & \dot f(t)=\gamma (u(t)(\Fmax-f(t)) -f(t)) \\
 & f(t_2)=\bar f, \qquad \tau \int_{t_2}^{t_f} f(t)^2 e^{-A(t-t_2)}\ dt=\frac{\sigma_f}A (1-e^{-A(t_f-t_2)})+
\gamma_1 e^0.\end{split}\eeq
 This problem can be kept as the full problem for the end of race. It provides a solution which is very close to that of Figure \ref{figfin}. Otherwise, one can try to reduce further the problem to have a simple expression for the velocity. In \cite{pess}, an approximation for such a problem by a sigmoid function is used. In our case, as computed in the Appendix, this yields the following sigmoid
\begin{equation}\label{endrace}
f(t) =\frac{ \Fmax}{ 1+  ( {\Fmax}/ {\bar f}-1) e^{-\gamma \lambda  \Fmax (t-t_2)}}
\end{equation}
where $\lambda$ is chosen such that the $L^2$ norm of $f$ satisfies  condition \eqref{eqenerA}. Then, since $v=\tau f$, this provides the final estimate for the velocity.
 This estimate yields an increasing velocity at the end of the race. It does not capture
  the short decrease at the very end of the race. But this changes very slightly the estimate on $t_f-t_2$ or on the sprint velocity at the end and is not meaningful for a runner, so we can safely ignore it for our approximations.

  Once we have this final approximation for the velocity, we have to match the length of the turnpike central phase so that the integral of $v$ is exactly $d$, which yields \eqref{t2t1}. This reduces very slightly the turnpike phase from 194.64 seconds
 to 193.81 seconds for our simulations.

\hfill

Our distance is made up of 3 parts: the turnpike distance  which is totally determined by $\gamma_1$ and $\gamma_2$ and the distance run in the initial and final parts. Of course, since the sum is prescribed, only one of the two is free. So for instance, in the final phase if we determine the duration of this final phase by some estimate like above, the initial phase has to match the total distance, but nevertheless is safely  estimated from the turnpike.

\section{Comparison with a real $1500$ m}

 The runners' oxygen uptake was recorded in \cite{hanon2008pacing} by means of a telemetric gas exchange system. This allowed to observe that the \vo\ reached a peak
 in around $450$\,m from start, with a significant decrease between $450$ and $550$ meters. Then the \vo\ remained constant for $800$ meters, before a decrease of $10\%$ at the
 end of the race.
 To match more precisely the \vo\ curve of \cite{hanon2008pacing}, we add an extra piece to the curve of $\sigma$,
  before the long mean value $\bar \sigma$: after the initial increase, there is a local maximum before decreasing to the constant turnpike value:
\begin{equation*}
\sigma (e)=\left\{\begin{array}{ll}
\displaystyle \bar \sigma \frac{e}{e^0\gamma_1}+\sigma_f \left(1-\frac{e}{e^0 \gamma_1}\right) & \displaystyle \quad\hbox{if}\quad \frac{e}{e^0}<\gamma_1\\[3mm]
\displaystyle \bar\sigma  & \displaystyle \quad\hbox{if}\quad \gamma_1\leq \frac{e}{e^0}\leq \gamma_+ \\[3mm]
\displaystyle \bar \sigma+0.8 \frac{e-\gamma_+e^0}{e^0-\gamma_2-\gamma_+e^0} & \displaystyle \quad\hbox{if}\quad \frac{e}{e^0}\geq \gamma_+ \quad\hbox{and}\quad e^0-e>\gamma_2 \\[3mm]
\displaystyle (\bar \sigma+0.8 -\sigma_r) \frac{e^0 -e}{\gamma_2}+\sigma_r & \displaystyle \quad\hbox{if}\quad e^0 -e<\gamma_2
\end{array}\right.
\end{equation*}
We take roughly the same parameters as before except for $\gamma_2=2000$ and $\gamma_+=1-\gamma_2/e^0-400/e^0$.
The others are $\sigma_r=6$,  $\sigma_f=20$,  $\bar\sigma=22$, $\gamma_1=0.15$,
 $F_{\mathrm{max}}=8$, $\tau=1.032$,  $e^0=4651$,  $\gamma=0.0025$, $v_0=1$.
\begin{figure}[ht]
\begin{centering}
\includegraphics[width=.9\textwidth]{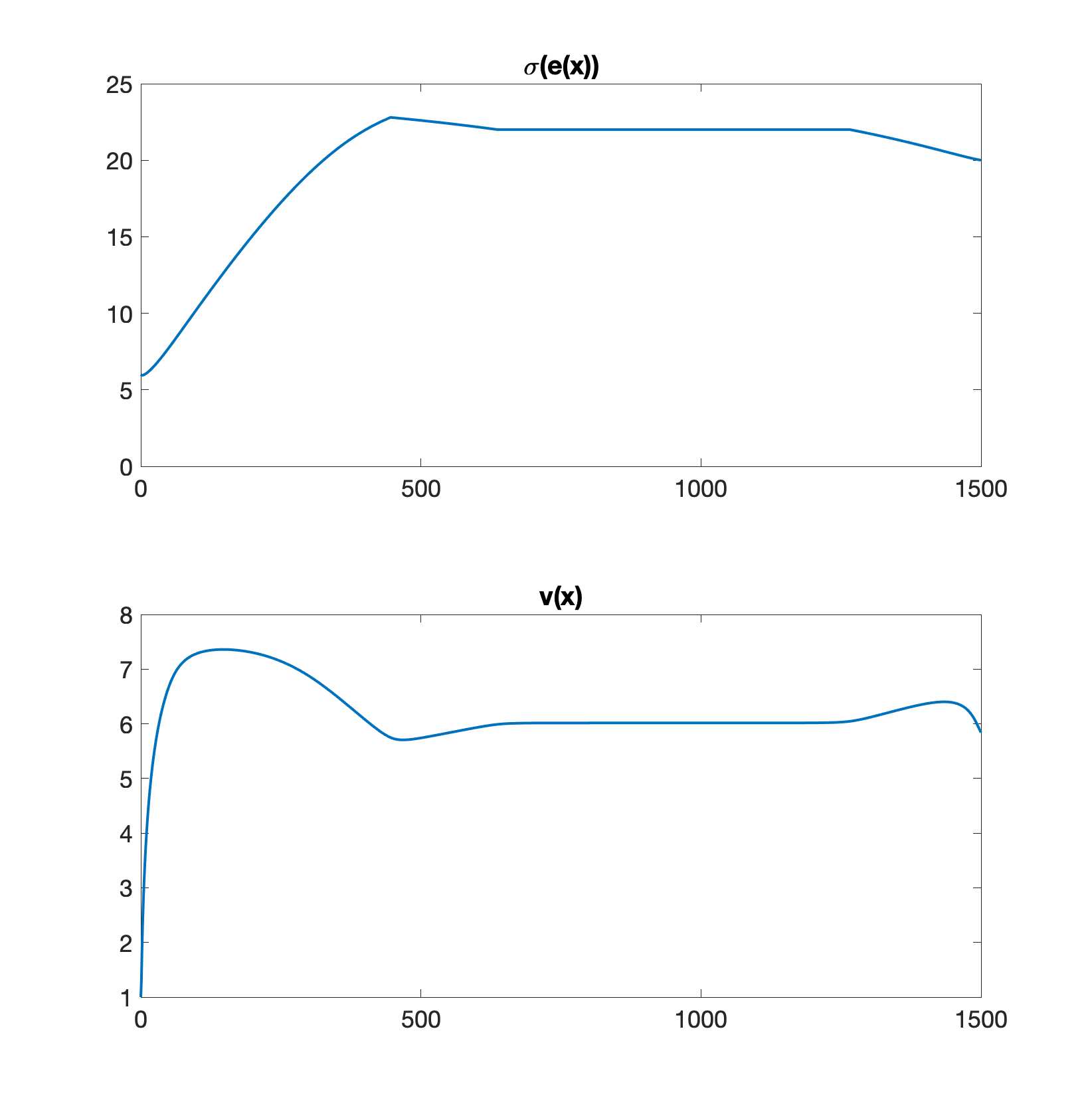}
\par\end{centering}
\caption{Modified $\sigma$ in four pieces and optimal velocity vs distance for a $1500$\,m.}\label{fighanon}
\end{figure}
 Then we see in Figure \ref{fighanon} that the velocity has a local minimum in the region where $\sigma$ has a local maximum,which matches exactly the velocity profile in \cite{hanon2008pacing}. Small variations in $\sigma$ always provide variations in the velocity profile with the opposite sense.

It is well known that successful athletes in a race are not so much those who speed up a lot at the end but those who avoid slowing down too much . We have noticed that if the maximal force at the beginning of the race is too high, then the velocity tends to fall down at the end of the race, leading to a bad performance. For a final in a world competition, it is observed in  \cite{hanley2019} that the best strategy is J-shaped, which means reaching maximal speed at the end of the race. But this is not available to all athletes.
 The runners profile of these simulations are not world champions but only successful in French regional races. Therefore, their pacing strategy is either  U-shaped (the start and the finish are quicker)  or reverse J-shaped (greater starting pace).  This is very dependent on the relative values of running economy $\tau$, anaerobic energy $e^0$ and profile of \vo.
 Moreover, top runners use pace variation according to laps as their winning tactics \cite{lapresa},
 but this is not active on the level of runners we have analyzed in this paper.







\section{Running uphill or downhill}
Our model also allows to deal with slope or ramps. Indeed, one has to change the Newton law of motion to take into account   a dependence on the slope $\beta(x)$ at distance $x$ from the start, which is the cosine of the angle. If we denote by $g$ the gravity,
the velocity equation changes into
$$
\dot v(t) = -\frac{v(t)}{\tau}+f(t)-g\beta(x(t)) .
$$
If the track goes uphill or downhill with a constant rate $\delta$, then  in the turnpike estimate, this becomes
$$
\bar v=\tau \bar f -g\tau \delta
$$
where $\delta$ is positive when the track goes up and negative when it goes down.
If the slope is constant for the whole race, the turnpike estimate can be computed.

If we assume a slope $\beta(x)$ which is constant equal to $\delta$, the new turnpike estimate is
$$\bar v=\frac{(e^0-dg\delta) \tau}{2d} +\sqrt{\bar\sigma\tau+\left ( \frac{(e^0-dg\delta)\tau}{2d}\right )^2}.$$If the slope is small, one can make an asymptotic expansion in terms of $\delta$ to find the difference in velocity
$$
\triangle v=-g\delta \tau \left (\frac 12+\frac 1{\sqrt{\frac{\bar \tau}4+\bar\sigma\left(\frac{d}{e^0}\right)^2}}\right) .
$$
But if the slope is constant for a small part of the race, then the variation of velocity cannot be computed locally because the whole mean velocity of the race is influenced by a local change of slope as we will see in the last part of the paper.

\begin{figure}[ht]
\centerline{\includegraphics[width=16cm]{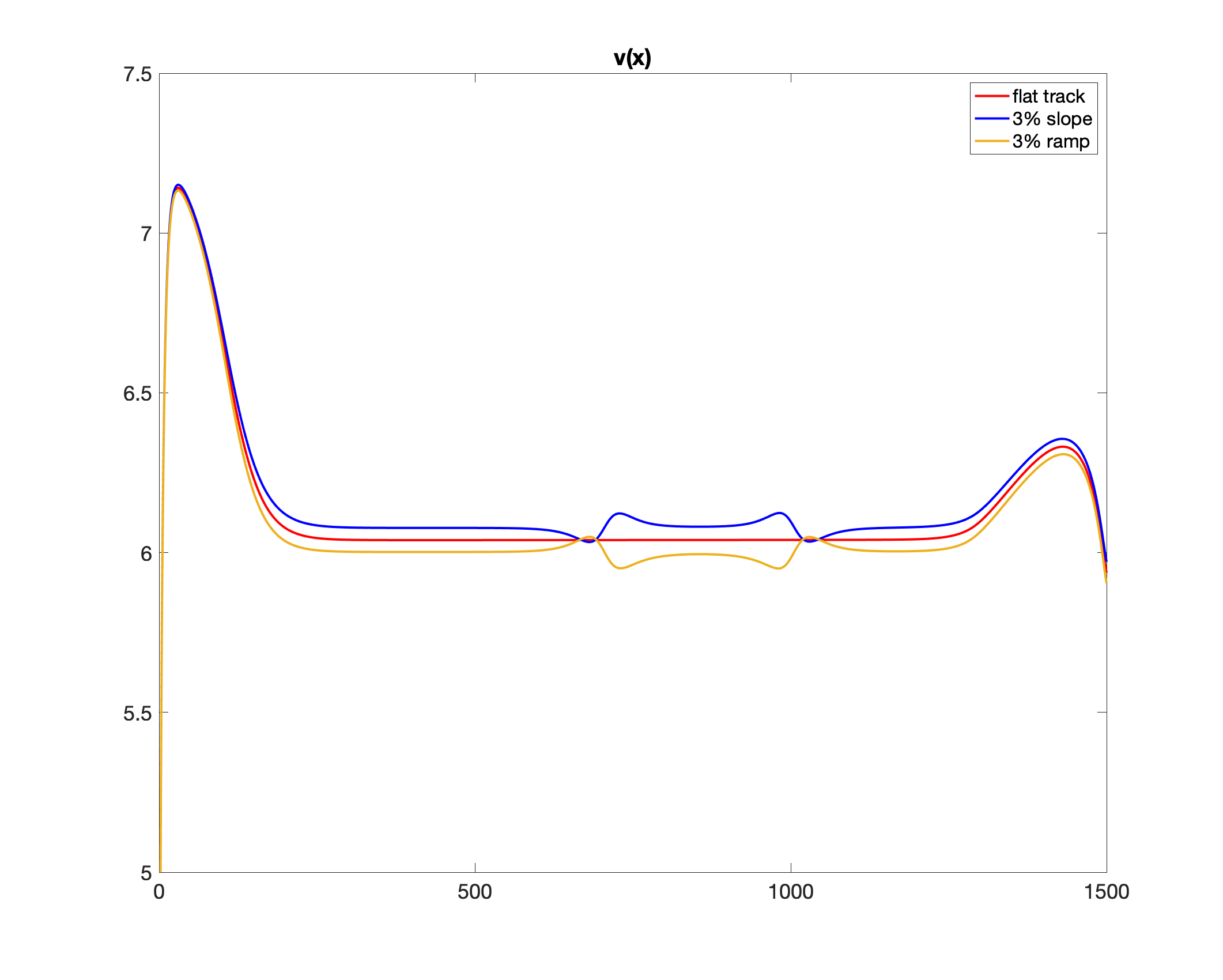}}
\caption{Velocity vs distance for a $1500$\,m, on a flat track (red), on a track with a $3\%$ slope between $700$\,m  and $1000$\,m (orange) and on a track with a  $3\%$ ramp between $700$\,m  and $1000$\,m  (blue).} \label{figcompp}
\end{figure}

Nevertheless, because the energy is involved, a change of slope, even locally implies a change of the turnpike velocity on the whole race.
We have chosen to put slopes and ramps of $3\%$ for $300$\,m. We see in Figure \ref{figcompp} that without slope we have an intermediate turnpike value, but with a slope or ramp even only for $300$\,m, the whole turnpike velocity is modified.

\begin{figure}[ht]
\centerline{\includegraphics[width=17.8cm]{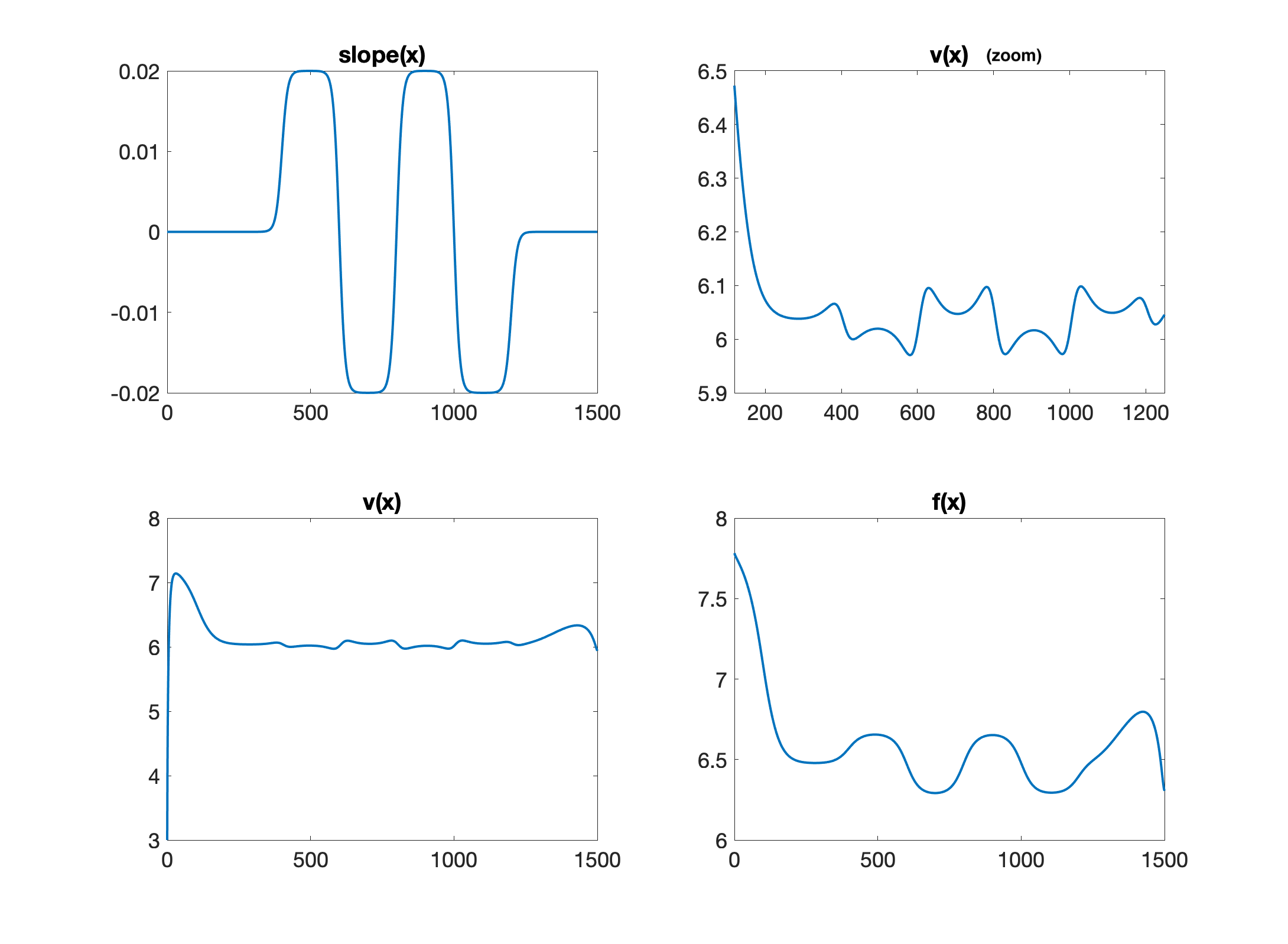}}
\caption{Slope, velocity, zoom on the velocity and force for a $1500$\,m with slopes and ramps. There is a slope  of $2\%$  between $400$\,m  and $600$\,m  and then between $800$\,m  and $1000$\,m. There is a ramp of  $2\%$  between $600$\,m  and $800$\,m  and then between $1000$\,m  and $1200$\,m.} \label{figpente}
\end{figure}

To illustrate further the slope effect, we have put a periodic slope and ramp of $200$\,m between $300$\,m and $1200$\,m. We use the same parameters as in the previous section. We see  in Figure \ref{figpente} that the turnpike velocity is affected.  When going down, a runner speeds at the end of the ramp, but his velocity has a local maximum at the middle of the ramp. Similarly, it has a local minimum at the middle of the slope. The variations in velocity are very small since they are of  order of a few percents. But this allows to understand that slopes and ramps are not local perturbations on the pacing profile.

\section{Conclusion}

We have provided a model for pace optimization. This involves a control problem in order to use the maximal available propulsive force and energy  to produce the optimal running strategy and minimize the time to run and the motor control. For sufficiently long races (above $1500$\,m), the optimal strategy is well approximated by a turnpike problem that we describe.   Simplified estimates for the peak velocity and velocity profiles related to aerobic, anaerobic energy and effect of the motor control are obtained and fit the simulations. The effect of the parameters and slope and ramps are analyzed. The potential applications of this turnpike theory would be to derive a simpler model for pacing strategy that could be encompassed in a running app. Indeed, the advantage of our simplified formulation for the velocity is that if we have velocity data of a runner on a race, and have access to his \vom, then we can infer the values of all the physiological parameters and therefore predict his optimal strategy on a fixed distance.

\section*{Appendix: Simplified motor control problem}
We want to study the simplified optimal control problem
\begin{equation*}
\begin{split}
& \min \int_0^T u (t)^2\, dt  \\
& \dot f(t) = \gamma (u(t)(\Fmax -f(t)) -f(t)) \qquad\qquad f(0)=\bar f\quad\textrm{and}\quad \int_0^T f(t)^2e^{-At}\ dt =\alpha ,
\end{split}
\end{equation*}
related to the one in \cite{pess} where there is no condition on the $L^2$ norm of $f$ but a final condition on $f(T)=F$ and a cost $\int_0^T u^2-k F$. In our case, we want to estimate $f(T)$ in terms of the parameters.

 The corresponding simplified problem for the beginning of the race is
 \begin{equation*}
\begin{split}
& \min \int_0^T u (t)^2\, dt  \\
& \dot f(t) = \gamma (u(t)(\Fmax -f(t)) -f(t)) \qquad\qquad f(T)=\bar f\quad\textrm{and}\quad \int_0^T f(t)^2e^{-At}\ dt =\alpha ,
\end{split}
\end{equation*}
 where we want to estimate $f(0)$ and understand why $f(t)$ is almost linear.
 Actually, at the beginning of the race the integral constraint would rather be of the form $\int_0^T f(t)v(t)\, dt=\alpha$ but this does not change the arguments developed hereafter.

 Because of the integral constraint on $f$, the above problem can be equivalently rewritten as
\begin{equation}\label{motou}
\begin{split}
& \min \int_0^T u (t)^2\, dt  \\
& \dot f(t) = \gamma (u(t)(\Fmax -f(t)) -f(t)) \qquad\qquad\qquad f(T)=\bar f , \\
& \dot y(t) = f(t)^2e^{-At} \qquad\qquad\qquad\qquad\qquad\qquad\qquad y(0)=0,\quad y(T)=\alpha .
\end{split}
\end{equation}
Let us apply the Pontryagin maximum principle to the optimal control problem \eqref{motou} (see \cite{LeeMarkus,Pontryagin,Trelat_book}).
Denoting by $p_f$ and $p_y$ the co-states associated, respectively, to the states $f$ and $y$, the Hamiltonian of the problem is
\begin{equation}\label{eqH}
H=p_f\gamma (u(\Fmax -f) -f) + p_y f^2e^{-At} -\frac{1}{2} u^2 .
\end{equation}
The condition $\frac{\partial H}{\partial u}=0$ yields $u=p_f\gamma (\Fmax -f)$.
Therefore,  the equation for $\dot f$ can be rewritten as
\begin{equation}\label{newdotf}
\dot f = \gamma \bigl (p_f \gamma (\Fmax -f)^2 -f \bigr ).
\end{equation}
 In order to estimate the solutions, we can assume that $p_f$ is not far from a constant
which allows an explicit integration of \eqref{newdotf}. Indeed the equation $p_f \gamma (\Fmax -f)^2 -f=0$ has two roots $f_1$ and $f_2$ and the solution of \eqref{newdotf} is thus the \emph{sigmoid} function
\begin{equation}\label{ftapp}
f(t)=f_2+\frac {f_1 -f_2}{1-\frac {\bar f -f_1}{\bar f -f_2} e^{\mu (t-T)}}
\end{equation}
with $\mu=p_f\gamma^2(f_1-f_2)$.
This allows to compute $f(0)$. Furthermore, if one approximates $e^{\mu (t-T)}$ by $1+\mu (t-T)$, then
$$
f(t)\simeq \bar f\frac{(\bar f -f_2)(\bar f -f_1)}{f_1-f_2}\mu (t-T)
$$
which is the linear approximation we have made for the first part of the race.

For the end of the race, the problem is similar except that it is an initial condition $f(0)=\bar f$ and we look for a final estimate on $f(T)$. A similar computation leads to the equivalent of \eqref{ftapp} which is the sigmoid function
\begin{equation}\label{ftappend}
f(t)=f_2+\frac {f_1 -f_2}{1-\frac {\bar f -f_1}{\bar f -f_2} e^{\mu t}},
\end{equation} which can also be rewritten as \eqref{endrace}.

\bibliographystyle{spmpsci}
\bibliography{biblio,manu}

\begin{thebibliography}{10}
\providecommand{\url}[1]{{#1}}
\providecommand{\urlprefix}{URL }
\expandafter\ifx\csname urlstyle\endcsname\relax
  \providecommand{\doi}[1]{DOI~\discretionary{}{}{}#1}\else
  \providecommand{\doi}{DOI~\discretionary{}{}{}\begingroup
  \urlstyle{rm}\Url}\fi

\bibitem{aft}
Aftalion, A.: How to run 100 meters.
\newblock SIAM Journal on Applied Mathematics \textbf{77}(4), 1320--1334 (2017)

\bibitem{AB}
Aftalion, A., Bonnans, J.F.: Optimization of running strategies based on
  anaerobic energy and variations of velocity.
\newblock {S}{I}{A}{M} {J}ournal on {A}pplied {M}athematics \textbf{74}(5),
  1615--1636 (2014)

\bibitem{AM}
Aftalion, A., Martinon, P.: Optimizing running a race on a curved track.
\newblock PloS one \textbf{14}(9), 0221572 (2019)

\bibitem{AT_RSOS}
Aftalion, A., Tr\'elat, E.: How to build a new athletic track to break records.
\newblock R. Soc. Open Sci. \textbf{200007}(April), 10 pp. (2020)

\bibitem{lapresa}
Arag{\'o}n, S., Lapresa, D., Arana, J., Anguera, M.T., Garz{\'o}n, B.: Tactical
  behaviour of winning athletes in major championship 1500-m and 5000-m track
  finals.
\newblock European Journal of Sport Science \textbf{16}(3), 279--286 (2016)

\bibitem{laumond}
Arechavaleta, G., Laumond, J.P., Hicheur, H., Berthoz, A.: An optimality
  principle governing human walking.
\newblock IEEE Transactions on Robotics \textbf{24}(1), 5--14 (2008)

\bibitem{behncke1993mathematical}
Behncke, H.: A mathematical model for the force and energetics in competitive
  running.
\newblock {J}ournal of {M}athematical {B}iology \textbf{31}(8), 853--878 (1993)

\bibitem{bsmall}
Behncke, H.: Small effects in running.
\newblock {J}ournal of {A}pplied {B}iomechanics \textbf{10}(3), 270--290 (1994)

\bibitem{BHKM}
Billat, V., Hamard, L., Koralsztein, J., Morton, R.: Differential modeling of
  anaerobic and aerobic metabolism in the 800-m and 1,500-m run.
\newblock J Appl Physiol. \textbf{107}(2), 478--87 (2009)

\bibitem{capo}
Bravo, M.J., Caponigro, M., Leibowitz, E., Piccoli, B.: Keep right or left,
  towards a cognitive-mathematical model for pedestrians.
\newblock Networks and Heterogeneous Media \textbf{10}, 559 (2015)

\bibitem{CasHan}
Casado, A., Hanley, B., Jimenez-Reyes, P., Renfree, A.: Pacing profiles and
  tactical behaviors of elite runners.
\newblock Journal of Sport and Health Science  (2020).
\newblock \doi{10.1016/j.jshs.2020.06.011}

\bibitem{fosterbeating}
Foster, C., de~Koning, J.J., Thiel, C., Versteeg, B., Boullosa, D.A., Bok, D.,
  Porcari, J.P.: Beating yourself: How do runners improve their own records?
\newblock International Journal of Sports Physiology and Performance pp. 1--10
  (2019)

\bibitem{Fourer2002}
Fourer, R., Gay, D.M., Kernighan, B.W.: AMPL: A mathematical programming
  language.
\newblock AT \& T Bell Laboratories Murray Hill, NJ 07974 (1987)

\bibitem{hh}
Hanley, B., Hettinga, F.J.: Champions are racers, not pacers: an analysis of
  qualification patterns of olympic and iaaf world championship middle distance
  runners.
\newblock Journal of sports sciences \textbf{36}(22), 2614--2620 (2018)

\bibitem{hanley2019}
Hanley, B., Stellingwerff, T., Hettinga, F.J.: Successful pacing profiles of
  olympic and iaaf world championship middle-distance runners across qualifying
  rounds and finals.
\newblock International journal of sports physiology and performance
  \textbf{14}(7), 894--901 (2019)

\bibitem{hanon2008pacing}
Hanon, C., Leveque, J.M., Thomas, C., Vivier, L.: Pacing strategy and {V}{O}2
  kinetics during a 1500-m race.
\newblock {I}nternational {J}ournal of {S}ports {M}edicine \textbf{29}(3),
  206--211 (2008)

\bibitem{hanon2011effects}
Hanon, C., Thomas, C.: Effects of optimal pacing strategies for 400-, 800-, and
  1500-m races on the {V}{O}2 response.
\newblock {J}ournal of {S}ports {S}ciences \textbf{29}(9), 905--912 (2011)

\bibitem{hettingabrian}
Hettinga, F.J., Edwards, A.M., Hanley, B.: The science behind competition and
  winning in athletics: using world-level competition data to explore pacing
  and tactics.
\newblock Frontiers in Sports and Active Living \textbf{1}, 11 (2019)

\bibitem{keller1974optimal}
Keller, J.B.: Optimal velocity in a race.
\newblock American {M}athematical {M}onthly pp. 474--480 (1974)

\bibitem{pess}
Le~Bouc, R., Rigoux, L., Schmidt, L., Degos, B., Welter, M.L., Vidailhet, M.,
  Daunizeau, J., Pessiglione, M.: Computational dissection of dopamine motor
  and motivational functions in humans.
\newblock Journal of Neuroscience \textbf{36}(25), 6623--6633 (2016)

\bibitem{LeeMarkus}
Lee, E.B., Markus, L.: Foundations of optimal control theory.
\newblock John Wiley \& Sons, Inc., New York-London-Sydney (1967)

\bibitem{mathis1989effect}
Mathis, F.: The effect of fatigue on running strategies.
\newblock {S}{I}{A}{M} {R}eview \textbf{31}(2), 306--309 (1989).
\newblock \urlprefix\url{http://www.jstor.org/stable/2030430 .}

\bibitem{AMH}
Mercier, Q., Aftalion, A., Hanley, B.: A model for world-class 10,000 m running
  performances: Strategy and optimization.
\newblock Frontiers in Sports and Active Living \textbf{2}, 226 (2021).
\newblock \doi{10.3389/fspor.2020.636428}.
\newblock
  \urlprefix\url{https://www.frontiersin.org/article/10.3389/fspor.2020.636428}

\bibitem{Per}
Peronnet, F., Massicote, D.: Table of nonprotein respiratory quotient: an
  update.
\newblock Can J Sport Sci \textbf{9}, 16--23 (1991)

\bibitem{Pontryagin}
Pontryagin, L.S., Boltyanskii, V.G., Gamkrelidze, R.V., Mishchenko, E.F.: The
  mathematical theory of optimal processes.
\newblock Translated from the Russian by K. N. Trirogoff; edited by L. W.
  Neustadt. Interscience Publishers John Wiley \& Sons, Inc.\, New York-London
  (1962)

\bibitem{thiel2012pacing}
Thiel, C., Foster, C., Banzer, W., De~Koning, J.: Pacing in olympic track
  races: competitive tactics versus best performance strategy.
\newblock Journal of sports sciences \textbf{30}(11), 1107--1115 (2012)

\bibitem{TJ}
Todorov, E., Jordan, M.I.: Optimal feedback control as a theory of motor
  coordination.
\newblock Nature Neuroscience \textbf{5}(11), 1226 (2002)

\bibitem{Trelat_book}
Tr\'{e}lat, E.: Contr\^{o}le optimal.
\newblock Math\'{e}matiques Concr\`etes. [Concrete Mathematics]. Vuibert, Paris
  (2005).
\newblock Th\'{e}orie \& applications. [Theory and applications]

\bibitem{trelat2020}
Tr{\'e}lat, E.: Linear turnpike theorem.
\newblock Preprint arXiv:2010.13605  (2020)

\bibitem{TZ}
Tr{\'e}lat, E., Zuazua, E.: The turnpike property in finite-dimensional
  nonlinear optimal control.
\newblock Journal of Differential Equations \textbf{258}(1), 81--114 (2015)

\bibitem{tucker2006non}
Tucker, R., Bester, A., Lambert, E.V., Noakes, T.D., Vaughan, C.L., Gibson,
  A.S.C.: Non-random fluctuations in power output during self-paced exercise.
\newblock {B}ritish {J}ournal of {S}ports {M}edicine \textbf{40}(11), 912--917
  (2006)

\bibitem{tucker2009physiological}
Tucker, R., Noakes, T.D.: The physiological regulation of pacing strategy
  during exercise: a critical review.
\newblock British Journal of Sports Medicine \textbf{43}(6), e1--e1 (2009)

\bibitem{Waechter2006}
W{\"a}chter, A., Biegler, L.T.: On the implementation of an interior-point
  filter line-search algorithm for large-scale nonlinear programming.
\newblock Mathematical Programming \textbf{106}(1), 25--57 (2006)

\end{thebibliography}

\end{document}